\documentclass[11pt, a4paper,twoside]{article}
\usepackage{amssymb,amsfonts,amsmath,amsthm,euscript}
\usepackage[english]{babel}
\usepackage[bf,small]{caption}
\usepackage{color}
\usepackage{dsfont}
\usepackage{geometry}
\usepackage{graphicx}
\usepackage{mathrsfs}
\usepackage{mathtools}
\usepackage{multicol}
\usepackage{multirow}
\usepackage{placeins}
\usepackage{setspace}
\usepackage{sectsty}
\usepackage{subfigure}
\usepackage{mdwlist}
\sectionfont{\large}
\allowdisplaybreaks 
\DeclareMathOperator*{\cov}{\mathrm{Cov}}

\usepackage[latin1]{inputenc}

\begin{document}
\theoremstyle{plain}
\newtheorem{thm}{Theorem}[section]
\newtheorem{cor}{Corollary}[section]
\newtheorem{prop}{Proposition}[section]
\newtheorem{lema}{Lemma}[section]

\theoremstyle{definition}
\newtheorem{df}{Definition}[section]
\newtheorem{rmk}{Remark}[section]
\newtheorem{exe}{Example}[section]

\addtolength{\parskip}{.2cm}
\numberwithin{equation}{section}
\numberwithin{table}{section}
\numberwithin{figure}{section}
\newcommand{\R}{\mathds{R}}
\renewcommand{\L}{\mathscr{L}}
\newcommand{\N}{\mathds{N}}
\newcommand{\Ns}{\mathds{N}^\ast}
\renewcommand{\t}{\theta}
\renewcommand{\d}{\partial}
\renewcommand{\qed}{{\hfill\footnotesize$\blacksquare$}\normalsize\\}
\newcommand{\p}{\varphi}
\renewcommand{\a}{\alpha}
\renewcommand{\l}{\lambda}
\newcommand{\bs}[1]{\boldsymbol{#1}}
\newcommand{\gm}{\gamma}
\newcommand{\tr}{\mathsf{T}}
\pagestyle{myheadings} 
\markboth{Parameterization of copulas and covariance decay }{G. Pumi and S.R.C. Lopes} 
\thispagestyle{empty}
{\centering
\huge{\bf  Parameterization of Copulas and Covariance Decay of Stochastic Processes}\vspace{.8cm}\\
\large{ {\bf Guilherme Pumi and S\'ilvia R.C. Lopes } \vspace{.1cm}\\
Mathematics Institute \\
Federal University of Rio Grande do Sul\vspace{.3cm}\\
This version: 05/02/2022\\
}
}
\vskip.6cm

\begin{abstract}
In this work we study the problem of constructing stochastic processes with a predetermined covariance decay by parameterizing its marginals and a given family of copulas. We show that the proposed methodology is compatibility-free and present several examples to illustrate the theory, including the important Gaussian and Euclidean families of copulas. We associate the theory to common applied time series models.
\vspace{.2cm}
\end{abstract}

\section{Introduction}

Let $\{X_t\}_{t=0}^\infty$ be a weakly stationary process. Within a probabilistic point of view, associated to any process, there is a sequence of distributions $\{F_n\}_{n=0}^\infty$, and to every $n$-upple of random variables $(X_{t_1},\cdots,X_{t_n})$ from the process, there is a copula $C_{t_1,\cdots,t_n}$ (Nelsen, 2006) associated to it. While the marginal distribution contains all the probabilistic information about the random variable itself, the copula account for the dependence among the variables $(X_{t_1},\cdots,X_{t_n})$. In a weakly stationary process, the covariance between $X_t$ and $X_{t+h}$ depends only on the lag $h$ and according to the copula version of  Hoeffding's lemma, the knowledge of the covariance between $X_t$ and $X_{t+h}$ can be obtained from the copula and the marginals of $X_t$ and $X_{t+h}$. Hence, studying the distributional properties among the variables is of vital importance in understanding the process' dependence structure.

Many classes of models are known to present certain covariance decay to zero as the lag $h$ increases. In the class of ARFIMA$(p,d,q)$  processes (Taqqu, 2003), for instance, the covariance may decay exponentially fast to zero or very slowly, depending on the parameter $d$. For heteroskedastic models such as ARCH (Engle, 1982) and GARCH (Bollerslev, 1986), the process itself is uncorrelated, but the absolute value and/or the square of the process are not and often present slow decay of covariance. Estimating the rate of such decay is generally of importance as it contains essential information about the long run structure of the process.

In recent years, the study of the dependence on univariate time series by means of copula has received growing attention. Darsow et al. (1992) provide the grounds of the modern use of copulas in stationary Markov processes. Some methods for constructing short memory time series based on conditional copulas are discussed in chapters 4 and 8 of Joe (1997).  The author also discusses methods for constructing Markovian processes and short memory time series based on parameterization of distributions belonging to a convolution-closed infinitely divisible class. Recent literature related to the present work includes Chen \& Fan (2006) and Chen et al. (2009) where semiparametric estimation in copula-based one dimensional stationary Markovian process are discussed. Lager{\aa}s (2010) discusses some non-standard behavior presented by some copula-based Markov processes. Beare and Seo (2012) consider the property of time reversibility in the context of copula-based Markov processes. Ibragimov (2009) studies higher order Markov processes in terms of copulas and study conditions under which a copula-based Markov process of some given order exhibit the so-called $m$-dependence, $r$-independence and conditional symmetry. Ibragimov \& Lentzas (2009) discuss a non-standard definition of long-range dependence in terms of copulas by considering the slow decay of some copula-based dependence measure other than the autocovariance/autocorrelation. Other interesting properties of copula-based Markov chain such as geometrical ergodicity, $\rho$-mixing, $\beta$-mixing, among others, are discussed in Chen et al. (2009) and Beare (2010a,b) (see also references therein). Although copulas are mainly applied to model nonlinear dependence, the linear case of covariance decay is not completely understood. It is worth noticing that recent literature on the subject either considers only the case of copula-based Markov processes or focus on non-standard definition of long-range dependence.

In the present note we intend to fill in this gap by proving two theorems. The first one (\ref{mp}) connects the covariance decay in a univariate time series and parametric bivariate copulas associated to lagged variables. For a process $\{X_n\}_{n=0}^\infty$ with absolutely continuous marginal distributions $\{F_n\}_{n=0}^\infty$ and for a given parametric family of copulas $\{C_\t\}_{\t\in\Theta}$ satisfying some minor regularity conditions,  we investigate how to obtain a given decay of covariance by choosing a parameterization $\{\t_n\}_{n=1}^\infty$ and by assuming that the copulas related to the pair $(X_{n_0},X_{n_0+h})$ depends only on the lag $h$ and are given by $C_{\t_h}$. We are especially interested in the case of slow decay of covariance, typical of long-range dependent processes, but the theory presented here is much more general than that and covers any arbitrary decay of covariance as well as for both, stationary and non-stationary processes.  We show that under suitable simple conditions, the parameterization on the copula family will ultimately determine the covariance decay of the pair $(X_{n_0},X_{n_0+h})$ as $h$ increases. The other theorem (\ref{compat}) we shall prove show that our approach is totally free of the so-called compatibility problem. Besides being fundamental to our approach, its proof is also an interesting application of index sets theory. To illustrate the theory, we present examples including several families of copulas widely used in practice. Special attention is given to Gaussian processes.

The paper is organized as follows. In the next section we recall a few concepts and results use\-ful in the paper as well as we introduce the idea of the paper by a simple example and show that the proposed approach is compatibility-free. In Sec\-tions 3 and 4, we consider the special cases of the Archimedean family of copulas and the Extreme Value copulas, respectively. In Section 5 we consider a general theory including arbitrary decay of covariance and present a discussion of the important case of Gaussian processes. Conclusions and final remarks are presented in Section 6.

\section{Preliminaries}

In this section we recall a few concepts and results we shall need in what follows. An $n$-dimensional copula is a distribution function defined in the $n$-dimensional hypercube $I^n$, where $I\vcentcolon=[0,1]$, and whose  marginals are uniformly distributed.  The interest in copulas has grown enormously in the last 15 years or so, especially because of its flexibility in applications. Copulas have been successfully applied and widely spread in several areas. In finances, copulas have been applied in major topics such as asset pricing, risk management and credit risk analysis among many others (see the books by Cherubini et al., 2004 and McNeil et al., 2005 for details). In  hydrology the modeling of rainfalls and storms often employ copulas to describe joint features between variables in the models (Palynchuk and Guo, 2011 and references therein). In econometrics, copulas have been widely employed in constructing multidimensional extensions of complex models (see Lee and Long, 2009 and references therein). In statistics, copulas have been applied in all sort of problems, such as development of dependence measures, modeling, testing, just to cite a few (see Mari and  Kotz, 2001, Nelsen, 2006 and references therein). Curiously, the main result in the theory, the celebrated Sklar's theorem recalled below, dates back to the late nineteen fifties. See Durante and Sempi (2015) for a proof.
\begin{thm}[Sklar's Theorem]Let $X_1,\cdots,X_n$ be random variables with joint distribution function $H$ and mar\-gi\-nals $F_1, \cdots , F_n$, respectively. Then, there exists a copula $C$ such that,
\begin{equation*}
H(x_1,\cdots,x_n)=C\big(F_1(x_1),\cdots,F_n(x_n)\big), \ \ \ \mbox{for all}\ \  (x_1,\cdots,x_n) \in \R^n.
\end{equation*}
If the $F_i$'s are continuous, then $C$ is unique. Otherwise, $C$ is uniquely determined on
$\mathrm{Ran}(F_1)\times\cdots\times\mathrm{Ran}(F_n)$. The converse also holds. Furthermore,
\begin{equation*}C(u_1,\cdots,u_n)=H\big(F_1^{(-1)}(u_1),\cdots,F_n^{(-1)}(u_n)\big),\ \ \  \mbox{for all} \ \  (u_1,\cdots,u_n)\in I^n,\end{equation*}
where for a function $F$, $F^{(-1)}$ denotes its pseudo-inverse given by $F^{(-1)}(x)\vcentcolon=\inf\big\{u\in \mathrm{Ran}(F): F(u)\geq x\big\}$ and, for a function $f$, $\mathrm{Ran}(f)$ denotes the range of $f$.
\end{thm}
The usefulness of Sklar's Theorem (and of copulas for extension) as a tool for statistical analysis and modeling lies on allowing one to deal separately with the joint dependence structure, characterized by the copula, and the marginals of a given random vector. This flexibility has been extensively explored in the literature. Another important result which shall be frequently used here is the copula version of the Hoeffding's lemma.
\begin{lema}[Hoeffding's Lemma]\label{hlema}
Let $X$ and $Y$ be two continuous random variables with marginal distributions $F$ and $G$, respectively, and copula $C$. Then, the covariance between $X$ and $Y$ is given by
\begin{equation}\label{hoef}
\cov(X,Y)=\iint_{I^2}\frac{C(u,v)-uv}{F'\big(F^{(-1)}(u)\big)G'\big(G^{(-1)}(v)\big)}\,dudv.
\end{equation}
\end{lema}
More details on the theory of copulas can be found in Joe (1997), Nelsen (2006) and Durante and Sempi (2015).

In this work,  $\N$ denotes the set of natural numbers  $\N\vcentcolon=\{0,1,2,\cdots\}$,  while $\N^\ast\vcentcolon=\N\setminus\{0\}$. For a given set $A\subseteq\R$, $\overline A$ denotes the closure of $A$ and $A'$ denotes the set of all accumulation points. For a vector $x\in\R^k$, $x^\tr$ denotes the transpose of $x$. The measure space behind the notion of measurable sets and functions is always  assumed (without further mention) to be $\big(\R^n,\mathscr B(\R^n),\mathfrak m\big)$ (or some appropriate restriction of it), where $\mathscr B(\R^n)$ denotes the Borel $\sigma$-field in $\R^n$ and $\mathfrak m$ is the Lebesgue measure in $\R^n$. Recall that a function $L:S\rightarrow\R$, for $S\subseteq\R$, is called \emph{slowly varying at} $a\in S'$ if $L$ is measurable, limited on a bounded interval and satisfies $\lim_{x\rightarrow a}L(cx)/L(x)=1$, for all $c>0$.  The set of slowly varying functions on the infinity will be denoted by
\begin{equation*}\L\vcentcolon=\{L: L \mbox{ is measurable, limited on a bounded interval and} \lim_{x\rightarrow \infty}L(cx)/L(x)=1,\, \forall c>0\}.\end{equation*}
We notice that if $L_1,L_2\in\L$, $L_1+L_2$, $L_1L_2\in\L$ and, for a constant $k\in\R\setminus\{0\}$, $kL_1\in\L$. Moreover, for any $\beta\in(0,1)$ and $L\in\L$, $L(n)n^{-\beta}\rightarrow 0$, fact usually explored in defining the concept of long-range dependence. More details on slowly varying function can be found in Bingham et al. (1987). For two functions $f$ and $g$, we denote $f(n)\sim g(n)$ if $f(n)/g(n)\rightarrow 1$, as $n$ goes to infinity. For an arbitrary sequence of real numbers $\{a_n\}_{n\in\Ns}$, consider the function $\psi_{a_{n}}:(0,\infty)\rightarrow\R$ given by $\psi_{a_{n}}(x)\vcentcolon=a_{\lceil x\rceil}$, where $\lceil \cdot \rceil$ denotes the ceiling function. For $L\in\L$, as a convention, by writing $a_nL\in\L$ we mean $\psi_{a_{n}}L\in\L$. We shall also write $a_n\in\L$ to mean $\psi_{a_n}\in\L$.  If $a_n\rightarrow a\neq0$, then always $a_n\in\L$. If $a_n\rightarrow 0$, then $a_n$ may or may not belong to $\L$, depending on the particular  convergence rate of the sequence to 0. For instance, $a_n=1/\ln(n)\in\L$, but $a_n=1/n\notin\L$.

As for long-range dependence, even though the literature on the matter is vast, there still no globally accepted definition for it. The most common definitions and their differences and similarities are discussed in section 4 of Taqqu (2003). In this work, we shall adopt the following general definition.
\begin{df}
We say that a weakly stationary process with finite variance $\{X_t\}_{t\in\N}$ presents \emph{long-range dependence} if
\begin{equation*}\cov(X_t,X_{t+h})\sim L(h)h^{-\beta}, \quad\mbox{ as }h\rightarrow\infty,\end{equation*}
for some $\beta\in(0,1)$ and $L\in\L$.
\end{df}
We start with an example to motivate the ideas we shall develop in the rest of the paper.
\begin{exe}\label{fgm}
Consider the \emph{Farlie-Gumbel-Morgenstern} (FGM, for short) family of copulas which consists of parametric copulas of the form
\begin{equation*}C_\t(u,v)=uv\big(1+\t(1-u)(1-v)\big),\end{equation*}
for $|\t|\leq 1$. As a particular case, we have the \emph{independence copula} $\Pi(u,v)\vcentcolon=uv$, for $\t=0$. Recall that a distribution function $F$ belongs to the so-called \emph{Type I Extreme Value} family with parameters $(a,b)\in\R\times(0,\infty)$, denoted by EVI$(a,b)$, if
\begin{equation*}F(x)=e^{-e^{-(x-a)/b}},\quad \mbox{ for all }x\in\R.\end{equation*}
Notice that $F'(F^{-1}(u))=b\hspace{1pt}u\ln(u)$, for all $u\in I$. Let $F_n\sim$ EVI$(a_n,b)$, for $\{a_n\}_{n\in\N}$ an arbitrary sequence of real numbers and $b>0$. Consider the FGM family of copulas with parameterization $\t_n\vcentcolon= \kappa_0^{-1}n^{-\alpha}$, $n\in\Ns$, for $\alpha>1$, where $\kappa_0\geq\zeta(\alpha)$ and $\zeta(\alpha)=\sum_{k=1}^\infty k^{-\alpha}$ is the so-called Riemman zeta function. For a fixed $n>2$, consider the collection of $\tbinom{n}{2}$ copulas  $C_{ij}\vcentcolon=C_{\t_{|i-j|}}$ for $i,j\in\{1,\cdots,n\}$, $i\neq j$.  Now consider the $n$-dimensional copula (cf. example 3.2 in Dolati and \'Ubeda-Flores, 2005)
\begin{equation}\label{alidolati}
C_n(u_1,\cdots,u_n)\vcentcolon=\sum_{1\leq i<j\leq n} \!\!\!\!\!C_{ij}(u_i,u_j)\!\prod_{\substack{k=1\\ k\neq i,j}}^n\!u_k-\frac{(n-2)(n+1)}{2}\prod_{l=1}^nu_l
\end{equation}
whose marginals are $C_{ij}$. Let $\{X_n\}_{n\in\N}$ be a sequence of random variables such that $X_n$ has distribution $F_n$, for each $n\in\N$.  For any $n>2$, let $C_n$ given in \eqref{alidolati} be the copula associated to $(X_0,\cdots, X_{n-1})$, which implies that the copula related to $(X_r,X_s)$ is $C_{rs}$, $r,s\in\{0,\cdots,n-1\}$ and $r\neq s$. Furthermore, Hoeffding's lemma  implies
\begin{equation*}\cov(X_t,X_{t+h})=\t_h\left(\int_I\frac{1-u}{b\ln(u)}\,du\right)^2=\frac{\ln(2)^2}{b^2}\,\t_h=\frac{\ln(2)^2}{b^2\kappa_0} h^{-\alpha},\end{equation*}
where the second equality follows from formula 4.267.8 in Gradshteyn and Ryzhik (2000). Since the construction \eqref{alidolati} is valid for any $n>2$, Sklar's theorem guarantees the existence of all finite dimensional distribution functions with the marginals, bivariate copulas and $n$-dimensional copulas as specified in the construction above. Therefore, by the Kolmogorov's existence theorem, we have just constructed a weakly stationary process $\{X_n\}_{n\in\N}$ by reparameterizing  a certain family of parametric copulas and its marginals.
\end{exe}
The objective of this paper is to derive conditions for which a certain decay of covariance can be achieved by simply parameterizing a certain family of copula and its marginal. We shall be interested in a more general framework than the one presented in Example \ref{fgm} in the sense that the theory covers both, weakly and strongly stationary processes as well as non-stationary ones.  Observe that it suffices to develop the theory for the sequence $\{(X_0,X_n)\}_{n=0}^{\infty}$, as the case $\{(X_{n_0},X_{n_0+n})\}_{n=0}^{\infty}$ for fixed $n_0\in\Ns$ is handled analogously by a simple re-indexation.

For a sequence $\{F_n\}_{n\in\N}$ of absolutely continuous distribution functions and for a sequence $\{C_n\}_{n\in\Ns}$ of copulas, consider  $\{X_n\}_{n\in\N}$ a sequence of random variables such that $X_n$ is distributed as $F_n$, for each $n\in\N$, and such that the copula associated to $(X_0,X_n)$ is $C_n$, for $n\in\Ns$. This construction is always possible by Sklar's theorem (notice that we are not making any joint distributional assumption, other than the ones implied to $(X_0,X_n)$, $n\in\Ns$). In this setting, if $\cov(X_0,X_n)\sim L(n)n^{-\beta}$, for $L\in\L$ and $\beta\in(0,1)$, we say that $(C_n,F_n)$ is compatible with long-range dependence structure. A strongly stationary long-range dependent process always satisfies this condition, but the construction alone is clearly not sufficient to specify a stochastic process, in which case we have to proceed analogously to Example \ref{fgm}. 

In principle, the sequence of random variables with the prescribed bidimensional copulas above need not to exist and one has to proceed carefully. This is the so-called compatibility problem, name given to the fact that, predefined a set of multivariate marginal distributions for a random vector $(X_1,\cdots,X_n)$, it may not exist an $n$-dimensional distribution with the given lower dimensional marginals. In this direction, we shall show that given a set of bidimensional copulas $C_1,\cdots,C_h$ corresponding to the random variables $(X_0,X_1),\cdots,(X_0,X_h)$, there \textbf{always} exists an $(h+1)$-dimensional copula associated to $(X_0,\cdots,X_h)$ with the given prescribed copulas as marginals. We observe that these can be viewed as a base construction of the so-called pair-copulas construction (Bedford and Cooke, 2001, 2002), which is based on pre-assignment of the conditional  bidimensional copula structure related to the $n$-dimensional copula under construction. The theorem presented here, however, has merit of its own first by solving the proposed problem; second, by doing so in a rather different fashion; and third, by not assuming any implicit conditional dependence construction as in the case of pair-copula construction. In order to do that, we shall show that our framework fits the context of Tiit (2002). First we need some concepts.
\begin{itemize*}
\item An \emph{index set} is a finite set of integers. For $I_1,\cdots, I_n$ index sets, the set $H=\{I_1,\cdots,I_n\}$ is called an \emph{index class}. $H$ is regarded as a set of sets.
\item The dimension of an index class $H$, denoted by $m(H)$, is the total number of distinct integers in the index sets contained in $H$.
\item The $s$-power of $H$, denoted by $s(H)$, is the number of index sets in $H$.
\item The $i$-power of $H$, denoted by $w(H)$, is the total number of integers contained in the index sets in $H$.
\item An index class is \emph{regular} if there is no $i\neq j\in\{1,\cdots,n\}$ such that $I_i\subset I_j$.
\item A regular index class $H$ is called \emph{divisible}, if there exists disjoint subsets $H_1$ and $H_2$ such that $H=H_1\cup H_2$ and for which there is no common integer in the index sets contained in $H_1$ and $H_2$. If these conditions fail, $H$ is called a \emph{regular indivisible index class}.
\item If in addition, the index sets in a regular indivisible index class $H$ satisfy $\max\{I_i\}\leq\max\{I_j\}$ whenever $i<j$, then the class is called a \emph{monotonic regular indivisible index class}, MRND for short.
\end{itemize*}
Tiit (2002) studies index sets in the following context. Suppose we want to construct an $n$-dimensional distribution for a random vector $(X_1,\cdots,X_n)$ given $q$ marginal distributions corresponding to vectors, say, $\{(X_{i_{k_1}},\cdots X_{i_{k_p}})\}_{p=1}^{q}$. In this situation, to each vector for which the distribution is given, it is attributed an index set with the indexes of the random variables corresponding to the given distribution. That is, the author considers the index class $H=\{I_1,\cdots,I_p\}$, where $I_j=\{i_{k_1},\cdots,i_{k_j}\}$. The next result is due to Tiit (2002) and is the key result in proving Theorem \ref{compat}. To simplify the notation, for a random vector $(X_1,\cdots,X_n)$ and an index set $I_q=\{i_1,\cdots,i_q\}\subseteq\{1,\cdots,n\}$, we shall denote by $P(I_j)$ the $p_j$-dimensional distribution of $(X_{i_1},\cdots,X_{i_{p_j}})$.
\begin{thm}\label{tiit} A sufficient condition for the existence of an $n$-dimensional distribution given a set of $p_j$-dimensional margins $\{P(I_j)\}_{j}$, where $H=\{I_j\}_j$ is a MRND class of index sets, is that
\begin{equation}\label{tiitcond}
w(H)=m(H)+s(H)-1,
\end{equation}
provided that the marginals corresponding to the indexes $I_i\cap I_k$, $i\neq k$, are the same.
\end{thm}
The proof of Theorem can be found in Tiit (2002). The next theorem assures the compatibility-free nature of the approach developed in the sections to come.
\begin{thm}\label{compat}
Let $C_0,\cdots,C_n$ be an arbitrary set of bidimensional copulas. Then there always exists an $n+1$-dimensional copula associated to a random vector $(X_0,\cdots, X_n)$ for which the copula associated to $(X_0,X_k)$ is $C_k$, $k\in{1,\cdots,n}$.
\end{thm}
\begin{proof}
We have $H=\big\{\{0,1\},\{0,2\},\cdots,\{0,n\}\big\}$.  It is clear that $H$ is an MRND class and since, in this case, $I_i\cap I_j=\{0\}$, $i\neq j$, the last condition on Theorem \ref{tiit} is immediately fulfilled. We have $w(H)=2n$, $m(H)=n+1$ and $s(H)=n$, so that \eqref{tiitcond} holds. The result now follows from Theorem \ref{tiit}.
\end{proof}
We emphasize that our approach is based on the copulas related to $\{(X_{0},X_h)\}_{h=1}^{\infty}$ only, so we do not make any assumption on the copulas related to any other pair of random variable, as this is not necessary to study the covariance decay of $\{(X_{0},X_h)\}_{h=1}^{\infty}$. The approach is actually tailored to avoid tackling the hard compatibility problem that may arise if one wants to have more control on the copulas other than $\{(X_{0},X_h)\}_{h=1}^{\infty}$. In this context we are only assuming minimal knowledge of the process' dependence structure and this minimality is exactly what allows the approach to work. We also make no conditional assumption as in the case of pair-copula construction.

To simplify the notation, for a sequence of absolutely continuous distribution  functions $\{F_n\}_{n\in\N}$ and a parametric copula $C_\t$, we shall write
\begin{equation*}l_n(x)\vcentcolon=F_n'\big(F_n^{(-1)}(x)\big)\qquad\text{and}\qquad \d^k C_\t(u,v)\vcentcolon=\frac{\d^k C_\t(u,v)}{\d\t^k}\,.\end{equation*}
\section{Archimedean Family}
Recall that the \emph{Archimedean family} of copulas consists of copulas of the form
\begin{equation*}C(u,v)=\p^{-1}(\p(u)+\p(v)),\end{equation*}
for some continuous decreasing convex possibly parametric function $\p:I\rightarrow[0,\infty]$,  $I\vcentcolon=[0,1]$, such that $\p(1)=0$, called the \emph{Archimedean generator} of the copula. We are interested in determining conditions in which a sequence of Archimedean copulas together with a sequence of absolutely continuous distribution functions is compatible with long-range dependence structure.
\begin{prop}\label{arch}
Let $\{\p_\t\}_{\t\in\Theta}$, for $\Theta\subseteq\R$ with non-empty interior, be a family of Ar\-chi\-me\-de\-an generators and $\{C_\t\}_{\t\in\Theta}$ be the correspondent Archimedean family. Suppose that there exists $a\in\Theta'$ such that $\lim_{\t\rightarrow a}\p_\t(t)=-\ln(t)$, where the limit is to be understood as the adequate lateral limit if $a\notin\mathrm{int}(\Theta)$. Also assume that there exists a set $D\subseteq\Theta$ with non-empty interior such that $a\in D'$ and $\p_\t$, seen as a function of $\t$, is of class $C^2$ in D. Let $\{F_n\}_{n\in\N}$ be a sequence of absolutely continuous distribution functions and define the sequences
\begin{equation*}K_1(n)\vcentcolon=\iint_{I^2}\lim_{\underset{\t\in D}{\t\rightarrow a}}\frac{\d C_{\t}(u,v)}{l_0(u)l_n(v)}dudv\quad\text{ and }\quad K_2(n)\vcentcolon=\iint_{I^2}\lim_{\underset{\t\in D}{\t\rightarrow a}}\frac{\d^2 C_{\t}(u,v)}{l_0(u)l_n(v)}dudv.\end{equation*}
Let $\{\t_n\}_{n\in\Ns}$ be a sequence in $D$ converging to $a$. Suppose that, for some  $\beta\in(0,1)$, $K_1(n)(\t_n-a)\sim L_1(n)n^{-\beta}$ and $K_2(n)(\t_n-a)^2=o(L_2(n)n^{-\beta})$, as $n$ goes to infinity, for $L_1,L_2\in\L$. Then, $(C_{\t},F_n)$ is compatible with long-range dependence structure.
\end{prop}
\begin{proof}
We present the proof for the case where $a\notin\mathrm{int}(\Theta)$ and assuming that $a>x$ for all $x\in D$. The other cases are proved analogously. Let $\{\a_m\}_{m\in\Ns}$ be a sequence of parameters in $D$ converging from the left to $a$. Applying a Taylor's expansion with Lagrange's remainder in $C_\t$ around $\t=a$, we obtain
\begin{align}\label{arch1}
C_{\t}(u,v)&=\lim_{m\rightarrow\infty}C_{\a_m}(u,v)+\lim_{m\rightarrow\infty}\d C_{\a_m}(u,v)(\t-a)+\frac12\lim_{m\rightarrow\infty} \d^2 C_{\a_m}(u,v)(\t_0-a)^2\nonumber\\
&= uv+\lim_{m\rightarrow\infty}\d C_{\a_m}(u,v)(\t-a)+\frac12\lim_{m\rightarrow\infty}\d^2 C_{\a_m}(u,v)(\t_0-a)^2,
\end{align}
where $\t_0\in[\t_n,a)$. Substituting $\t$ by $\t_n$ in \eqref{arch1}, we obtain
\begin{equation}\label{arch2}
C_{\t_n}(u,v)=uv+\lim_{m\rightarrow\infty}\d C_{\a_m}(u,v)(\t_n-a)+\frac12\lim_{m\rightarrow\infty}\d^2 C_{\a_m}(u,v)\big(\t_0(n)-a\big)^2,
\end{equation}
where $\t_0(n)\in[\t_n,a)$ depends on $n$ and satisfies $\lim_{n\rightarrow\infty}\t_0(n)=a$. Under the  notation of the enunciate, we obtain
\begin{align}\label{arch3}
  \iint_{I^2}\frac{C_{\t_n}(u,v)-uv}{l_0(u)l_n(v)}\,dudv&=\,\left(  \iint_{I^2} \lim_{m\rightarrow\infty}\frac{\d C_{\a_m}(u,v)}{l_0(u)l_n(v)}\,dudv\!\right)\!\!(\t_n-a)+\nonumber\\
&\hspace{1.5cm}+\,\frac12\left(  \iint_{I^2}\lim_{m\rightarrow\infty}\frac{\d^2 C_{\a_m}(u,v)}{l_0(u)l_n(v)}\,dudv\right)\!\big(\t_0(n)-a\big)^2\nonumber\\
&=K_1(n)(\t_n-a)+\frac12K_2(n)\big(\t_0(n)-a\big)^2.
\end{align}
Since $\t_0(n)\in[\t_n,a)$,
\begin{equation}\label{arch4}
\big|K_2(n)\big|\big(\t_0(n)-a\big)^2\leq\big|K_2(n)(\t_n-a\big)^2|=o(L_2(n)n^{-\beta}),
\end{equation}
by the hypothesis on $K_2$. The result now follows from \eqref{arch3} in view of \eqref{arch4} and the hypothesis on $K_1$.
\end{proof}
\begin{exe}\label{amh}
The Ali-Mikhail-Haq (AMH, for short) family of copulas is  Archimedean  with generator $\p_\t(t)=(1-\t)^{-1}\ln\big([1+\t(t-1)]t^{-1}\big)$, for $|\t|\leq 1$. The copulas of the AMH family have the form
\begin{equation*}C_\t(u,v)=\frac{uv}{1-\t(1-u)(1-v)}\,.\end{equation*}
As a particular case, we have $C_0=\Pi$. Let $\{F_n\}_{n\in\N}$ be a sequence of distribution functions such that $F_n\sim\mathrm{Exp}(\l_n)$, for $\{\l_n\}_{n\in\N}$ a sequence of positive real numbers. Simple calculations show that $l_n(x)=(1-x)/\l_n$ and
\begin{equation}\label{amh1}
K_1(n)=\l_0\l_n\left(\int_Iu\,du\right)^2=\frac{\l_0\l_n}{4}\quad \mbox{and}\quad K_2(n)=2\l_0\l_n\left(\int_Iu(1-u)du\right)^2=\frac{\l_0\l_n}{18}.
\end{equation}
From \eqref{amh1}, there are many ways to take $\l_n$ and $\t_n$ in order to obtain the desired decay. For instance, for $\beta\in(0,1)$ and $L\in\L$, by taking $\l_n\in\L$ and $\t_n\sim n^{-\beta}$, or $\l_n\sim n^{-\beta}$ and $\t_n\rightarrow 0$ but such that $\t_n\in\L$, then Proposition \ref{arch} applies and we conclude that $(C_\t,F_n)$ is compatible with long-range dependence structure. If $\l_n=\l_0>0$, for all $n$, as it is the case for a strongly stationary process, then by taking $\t_n\sim L(n)n^{-\beta}$, the same conclusion is obtained through Proposition \ref{arch}.
\end{exe}
\begin{exe}\label{GB}
The \emph{Gumbel-Barnett} family of copulas is Archimedean with generator $\p_\t(x)=\ln(1-\t\ln(x))$, for $\t\in(0,1]$ and copulas given by
\begin{equation*}C_\t(u,v)=uv\exp(-\t\ln(u)\ln(v)).\end{equation*}
As a limiting case we have $\lim_{\t\rightarrow 0^+}C_\t=\Pi$. Let $F_n$ be distributed as EVI$(a_n,b_n)$, $n\in\N$, for $\{a_n\}_{n\in\N}$ an arbitrary sequence of real numbers and $\{b_n\}_{n\in\N}$ a sequence of positive real numbers. A simple calculation shows that
\begin{equation*}K_1(n)=\frac{1}{b_0b_n}\left(\int_{I}du\right)^2=\frac{1}{b_0b_n}\quad\mbox{and}\quad K_2(n)=\frac{1}{b_0b_n}\left(\int_I\ln(u)du\right)^2=\frac{1}{b_0b_n}\,.\end{equation*}
Therefore, the situation is similar to Example  \ref{amh}. For instance, for $\beta\in(0,1)$ and $L\in\L$, taking $1/b_n\in\L$ and $\t_n\sim n^{-\beta}$, or $b_n\sim n^{\beta}$ and $\t_n\rightarrow 0$ but such that $\t_n\in\L$, then Proposition \ref{arch} applies and we conclude that $(C_\t,F_n)$ is compatible with long-range dependence structure. If $1/b_n=k_0>0$ for all $n$, as it is the case for a strongly stationary process, then taking $\t_n\sim L(n)n^{-\beta}$, the same conclusion is obtained through Proposition \ref{arch}.
\end{exe}
\begin{exe}\label{frexe}
Consider the \emph{Frank} family of copulas with Archimedean generator $\p_\t(x)=-\ln\big((e^{-\t t}-1)(e^{\t}-1)^{-1}\big)$, for $\t\in\R\setminus\{0\}$, with corresponding copula given by
\begin{equation*}C_\t(u,v)=-\frac1\t\ln\left(1+\frac{(e^{-\t u}-1)(e^{-\t v}-1)}{e^{-\t}-1}\right).\end{equation*}
This is one of the most applied Archimedean  copula in statistics. Details of its nice properties can be found in Nelsen (2006) and references therein.   As a particular case we have $C_0=\Pi$.  Let $F_n$ be distributed as EVI$(a_n,b_n)$, $n\in\N$, for $\{a_n\}_{n\in\N}$ an arbitrary sequence of real numbers  and $\{b_n\}_{n\in\N}$ a sequence of positive real numbers. Upon applying formula 4.267.8 in Gradshteyn and Ryzhik (2000), it is routine to show that
\begin{align}\label{frank}
K_1(n)=\frac{1}{2b_0b_n}\left(\int_I\frac{1-u}{\ln(u)}\, du\right)^2=\frac{\ln(2)^2}{2b_0b_n}\quad\mbox{ and}\quad K_2(n)=\iint_{I^2}\frac{p(u,v)}{b_0u\ln(u)b_nv\ln(v)}\,dudv=\frac{k_0}{b_0b_n},
\end{align}
where
\begin{equation*}p(u,v)\vcentcolon=-u^2v^3-u^3v^2+\frac{2u^3v^3}{3}-\frac{uv^2}2-\frac{u^2v}{2}+\frac{3u^2v^2}{2}+\frac{uv}{6}+\frac{u^3v}{3}+\frac{uv^3}{3}\end{equation*}
and $ k_0=2\ln(3)^2/3-2\ln(2)\ln(3)+3\ln(2)^2/2$.  In view of \eqref{frank}, for $\beta\in(0,1)$ and $L\in\L$, taking $1/b_n\in\L$ and $\t_n\sim n^{-\beta}$, or $b_n\sim n^{\beta}$ and $\t_n\rightarrow 0$ but such that $\t_n\in\L$, then by Proposition \ref{arch} we conclude that $(C_\t,F_n)$ is compatible with long-range dependence structure. If $1/b_n=k_0>0$ for all $n$, as it is the case for a strongly stationary process, then taking $\t_n\sim L(n)n^{-\beta}$, the same conclusion holds.
\end{exe}
\section{Extreme Value Copulas}

Recall that $C$ is an  \emph{Extreme Value Copula} (EVC, for short) if there exists a copula $C_0$ such that
\begin{equation*}C(u,v)=\lim_{n\rightarrow\infty}C^n_0(u^{\frac1n},v^{\frac1n}),\end{equation*}
for all $(u,v)\in I^2$. A method to ``parameterize'' the EVC family is devised in Pickands (1981) (see also Nelsen, 2006, p.97). The construction show that a copula $C$ belongs to the EVC family if it can be written as
\begin{equation*}C_A(u,v)=\exp\!\bigg(\!\ln(uv)A\!\left(\frac{\ln(u)}{\ln(uv)}\right)\!\!\!\bigg),\end{equation*}
for some possibly parametric convex function $A:I\rightarrow [0.5,1]$, called \emph{dependence function}, satisfying $A(0)=A(1)=1$ and $A(t)\in[\max\{t,1-t\},1]$, for all $t\in I$. As a particular case, when $A\equiv1$, $C_A=\Pi$. For a parametric dependence function $A_\t$, we shall denote the corresponding EVC family by $C_\t\vcentcolon=C_{A_\t}$. Under some conditions, the EVC family is also compatible with long-range dependence structure. To simplify the notation, let
\begin{equation*}\d^k A_\t(x)\vcentcolon=\frac{\d^k}{\d t^k}A_t(x)\bigg|_{t=\t},\quad \mathcal A_\t(u,v)\vcentcolon=A_\t\! \left(\frac{\ln(u)}{\ln(uv)}\right),\quad \mbox{and}\quad \partial \mathcal A_\t(u,v)\vcentcolon=\frac{\partial \mathcal A_\t(u,v)}{\partial\t},\end{equation*}
for $k=1,2$.  Proposition \ref{evc} below presents conditions similar to Proposition \ref{arch} to obtain the compatibility of a EVC family with long-range dependence structure.
\begin{prop}\label{evc}
Let $\{A_\t\}_{\t\in\Theta}$, for  $\Theta\subseteq \R$ with non-empty interior, be a family of dependence functions and let $\{C_\t\}_{\t\in\Theta}$ be the associated \emph{EVC} family. Suppose that there exists $a\in\Theta'$ such that $\lim_{\t\rightarrow a} A_\t\equiv 1$ where the limit is to be understood as the adequate lateral limit if $a\notin\mathrm{Int}(\Theta)$. Also assume that there exists a set $D\subseteq\Theta$ with non-empty interior such that $a\in D'$ and $\p_\t$, seen as a function of $\t$, is of class $C^2$ in D. Let $\{F_n\}_{n\in\N}$ be a sequence of absolutely continuous distribution functions and define the sequences
\begin{equation*}K_1(n)\vcentcolon=\iint_{I^2}\lim_{\underset{\t\in D}{\t\rightarrow a}}\frac{\d C_{\t}(u,v)}{l_0(u)l_n(v)}dudv\quad\text{ and }\quad K_2(n)\vcentcolon=\iint_{I^2}\lim_{\underset{\t\in D}{\t\rightarrow a}}\frac{\d^2 C_{\t}(u,v)}{l_0(u)l_n(v)}dudv.\end{equation*}
Let $\{\t_n\}_{n\in\Ns}$ be a sequence in $D$ converging to $a$. Suppose that, for some  $\beta\in(0,1)$, $K_1(n)(\t_n-a)\sim L_1(n)n^{-\beta}$ and $K_2(n)(\t_n-a)^2=o(L_2(n)n^{-\beta})$, as $n$ goes to infinity, for $L_1,L_2\in\L$. Then, $(C_{\t},F_n)$ is compatible with long-range dependence structure.
\end{prop}
\begin{proof}
The proof goes along the same lines as the proof of Proposition \ref{arch}.
\end{proof}
It is easy to see that equivalent expressions to $K_1$ and $K_2$ in Proposition \ref{evc} are as follows:
\begin{equation}\label{evc1}
K_1(n)=\iint_{I^2}\bigg[\lim_{\underset{\t\in D}{\t\rightarrow a}}\d \mathcal{A}_{\t}(u,v)\bigg]\frac{ uv\ln(uv)}{l_0(u)l_n(v)}\,dudv
\end{equation}
and
\begin{equation}\label{evc2}
K_2(n)=\iint_{I^2}\bigg(\ln(uv)\bigg[\lim_{\underset{\t\in D}{\t\rightarrow a}}\d\mathcal{A}_{\t}(u,v)\bigg]^2\!\!+\lim_{\underset{\t\in D}{\t\rightarrow a}}\d^2\mathcal{A}_{\t}(u,v)\bigg)\frac{uv\ln(uv)}{l_0(u)l_n(v)}\,dudv.
\end{equation}
The above expressions seem cumbersome at first glance, they are often simpler to calculate because the limits in the integrand usually result in simple expressions.
\begin{exe}
Consider the following dependence function and the associated EVC family
\begin{equation*}A_\t(t)=1-\t t(1-t) \quad \mbox{ and }\quad C_\t(u,v)=uv\,\mathrm{exp}\left(\frac{\t\ln(v)\big(\ln(v)-\ln(uv)\big)}{\ln(uv)}\right),\end{equation*}
for $\t\in I$. This is known as the generator of the Tawn mixed model (cf. Mari and Kotz, 2001, p.96). Notice that $A_0(t)\equiv1$, $\d A_0(t)=-t(1-t)$ and $\d^2 A_0(t)=0$. Let $F_n$ be distributed as  EVI$(a_n,b_n)$, $n\in\N$, for $\{a_n\}_{n\in\N}$ an arbitrary sequence of real numbers and $\{b_n\}_{n\in\N}$ a sequence of positive real numbers. By using \eqref{evc1} and \eqref{evc2},
\begin{equation*}K_1(n)=\iint_{I^2}\frac{\ln(v)\big(\ln(v)-\ln(uv)\big)}{b_0\hspace{1pt}b_n\ln(uv)^2}\,dudv=-\frac{1}{b_0\hspace{1pt}b_n}\iint_{I^2}\frac{\ln(u)\ln(v)}{\big(\ln(u)+\ln(v)\big)^2}\,dudv= -\frac{1}{6\hspace{1pt}b_0\hspace{1pt}b_n},\end{equation*}
and
\begin{equation*}K_2(n)=\iint_{I^2}\frac{\left(\ln(v)\big(\ln(v)-\ln(uv)\big)\right)^2}{b_0\hspace{1pt}b_n\ln(uv)^3}\,dudv=-\frac{1}{b_0\hspace{1pt}b_n}\iint_{I^2}\frac{\ln(u)\ln(v)}{\big(\ln(u)+\ln(v)\big)^3}\,dudv= -\frac{1}{15\hspace{1pt}b_0\hspace{1pt}b_n}\,,\end{equation*}
where the integrals in both expressions are calculated by changing variables to $u=e^{-t}$, by applying integral by parts and by using formula 3.353.1 in Gradshteyn and Ryzhik (2000). Therefore, for $\beta\in(0,1)$ and $L\in\L$, taking $1/b_n\in\L$ and $\t_n\sim n^{-\beta}$, or $b_n\sim n^{\beta}$ and $\t_n\rightarrow 0$ but such that $\t_n\in\L$, then Proposition \ref{evc} applies and we conclude that $(C_\t,F_n)$ is compatible with long-range dependence structure. If $1/b_n=k_0>0$ for all $n$, by taking $\t_n\sim L(n)n^{-\beta}$ the same conclusion holds.
\end{exe}
\section{General Theory}\label{gensec}
In the previous sections we studied two particular cases of a more general theory to be developed in this section. The focus of the previous section was on the compatibility with long-range dependence structure. However, as it shall become clear in Theorem \ref{gen} below, with a minimum effort we can  extend the theory to cover any arbitrary decay of covariance in a wide class of parametric families of absolutely continuous copulas satisfying minimal regularity conditions.
\begin{thm}\label{gen}
Let $\{C_\t\}_{\t\in\Theta}$ be a family of parametric copulas, for $\Theta\subseteq \R$ with non-empty interior. Suppose that there exists a point $a\in\Theta'$ such that $\lim_{\t\rightarrow a}C_\t=\Pi$, where the limit is to be understood as the adequate lateral limit if $a\notin\mathrm{int}(\Theta)$. Also assume that there exist a set $D\subseteq\Theta$ with  non-empty interior such that $a\in D'$ and $C_\t$, seen as a function of the parameter $\t$, is of class $C^2$ in $D$. Let $\{F_n\}_{n\in\N}$ be a sequence of absolutely continuous distribution functions and define the sequences
\begin{equation*}K_1(n)\vcentcolon=\iint_{I^2}\lim_{\underset{\t\in D}{\t\rightarrow a}}\frac{\d C_{\t}(u,v)}{l_0(u)l_n(v)}\,dudv\quad\text{ and }\quad K_2(n)\vcentcolon=\iint_{I^2}\lim_{\underset{\t\in D}{\t\rightarrow a}}\frac{\d^2 C_{\t}(u,v)}{l_0(u)l_n(v)}\,dudv.\end{equation*}
Let $\{\t_n\}_{n\in\Ns}$ be a sequence in $D$ converging to $a$ and let $\{X_n\}_{n\in\N}$ be a sequence of random variables such that $X_n\sim F_n$, $n\in\N$ and the copula associated with $(X_0,X_n)$ be $C_{\t_n}$. Given a measurable function $R:\R\rightarrow\R$ satisfying $\lim_{n\rightarrow\infty}R(n)=0$, a sufficient condition for $\cov(X_0,X_n)\sim R(n)$, as $n$ tends to infinity, is that $K_1(n)(\t_n-a)\sim R(n)$ and $K_2(n)(\t_n-a)^2=o(R(n))$.
\end{thm}
\begin{proof}
We present the proof for the case where $a\notin\mathrm{int}(\Theta)$ and assuming that $a>x$ for all $x\in D$. The other cases are dealt analogously.  Let $\{\a_m\}_{m\in\Ns}$ be an arbitrary sequence of parameters in $D$ converging from the left to $a$. Applying a Taylor's expansion with Lagrange's remainder in $C_\t$ around $\t=a$, and proceeding as in the proof of Proposition \ref{arch}, we obtain
\begin{equation}
C_{\t_n}(u,v)=uv+\lim_{m\rightarrow\infty}\d C_{\a_m}(u,v)(\t_n-a)+\frac12\lim_{m\rightarrow\infty}\d^2 C_{\a_m}(u,v)\big(\t_0(n)-a\big)^2,
\end{equation}
with the notation as in the enunciate. By Hoeffding's lemma it follows that
\begin{align*}
\cov(X_0,X_n) &= \iint_{I^2}\frac{C_{\t_n}(u,v)-uv}{l_0(u)l_n(v)}\,dudv=\,\left(  \iint_{I^2} \lim_{m\rightarrow\infty}\frac{\d C_{\a_m}(u,v)}{l_0(u)l_n(v)}\,dudv\!\right)\!\!(\t_n-a)+\nonumber\\
&\hspace{1.5cm}+\,\frac12\left(  \iint_{I^2}\lim_{m\rightarrow\infty}\frac{\d^2 C_{\a_m}(u,v)}{l_0(u)l_n(v)}\,dudv\right)\!\big(\t_0(n)-a\big)^2\nonumber\\
&=K_1(n)(\t_n-a)+\frac12K_2(n)\big(\t_0(n)-a\big)^2,
\end{align*}
where $\t_0(n)\in[\t_n,a)$ depends on $n$ and satisfies $\lim_{n\rightarrow\infty}\t_0(n)=a$, provided the integrals exist. By hypothesis
\begin{equation*}\big|K_2(n)\big|\big(\t_0(n)-a\big)^2\leq\big|K_2(n)(\t_n-a\big)^2|=o(R(n)),\end{equation*}
and the result follows from the hypothesis on $K_1(n)$.
\end{proof}
\begin{exe}[\emph{Covariance Decay on Gaussian Processes}]\label{gauss}
Let $\phi$, $\Phi$ and $\Phi^{-1}$ denote the density, the distribution function and the quantile function of a standard normal random variable, respectively. Also let $\Phi_\rho$ denote the distribution function of a bivariate normal distribution with mean $(0,0)^\tr$ and variance-covariance matrix given by $\Omega\vcentcolon=\bigl(\begin{smallmatrix}1&\rho\\ \rho&1\end{smallmatrix} \bigr)$. The so-called \emph{Gaussian family of copulas} comprehend the copulas given by
\begin{equation*}C_\rho(u,v)=\Phi_\rho\big(\Phi^{-1}(u),\Phi^{-1}(v)\big),\end{equation*}
for $\rho\in[-1,1]$. As a particular case we have $C_0=\Pi$.  Let $F_n=\Phi$, for all $n\in\N$. For simplicity, for $u\in I$, we shall denote $q_u\vcentcolon=\Phi^{-1}(u)$. The chain rule yields
\begin{equation*}\frac{\d C_\rho(u,v)}{\d \rho}=\frac{\rho C_\rho(u,v)}{2\pi(1-\rho)^{\frac32}}+\frac1{2\pi\sqrt{1-\rho}}\iint_{S}\bigg[\frac{xy}{1-\rho^2}-\frac{(x^2+y^2-2\rho xy)}{(1-\rho^2)^2}\bigg]
e^{-\frac{\rho(x^2+y^2-2\rho xy)}{2(1-\rho^2)}}dxdy,\end{equation*}
where  $S\vcentcolon=S(u,v)=(-\infty,q_u]\times(-\infty,q_v]$. Thus,
\begin{equation*}\lim_{\rho\rightarrow 0} \frac{\d C_\rho(u,v)}{\d \rho}=\frac{1}{2\pi}\left(\int_{-\infty}^{q_u}ye^{-\frac{y^2}{2}}dy\right)\left(\int_{-\infty}^{q_v}ye^{-\frac{y^2}{2}}dy\right)=\phi(q_u)\phi(q_v),\end{equation*}
where the last equality follows from integration by parts. Hence $K_1(n)=1$.
We move to calculate $K_2$. Elementary calculus yields
\begin{align*}
\lim_{\rho\rightarrow 0} \frac{\d^2 C_\rho(u,v)}{\d \rho^2}&=uv -\iint_Sx^2(1-y^2)\phi(x)\phi(y)dxdy\\
&=uv-\Big[uv-q_uq_v\phi(q_u)\phi(q_v)\Big]=q_uq_v\phi(q_u)\phi(q_v).
\end{align*}
Now, since $F_n\sim \mathcal{N}(0,1)$,
\begin{equation*}K_2(n)=\bigg(\int_Iq_udu\bigg)^2=\bigg(\int_\R y\phi(y)dy\bigg)^2=0,\end{equation*}
where the second equality follows by changing variables to $y=q_u$. Therefore, if $\{X_n\}_{n\in\N}$ is a sequence of standard normal random variables and the copula associated with $(X_0,X_n)$  is $C_{\rho_n}$, the calculations above show that the parameterization we choose for the sequence $\{\rho_n\}_{n\in\Ns}$ will ultimately determine the decay of covariance (and the decay of correlation for that matter, since $X_n\sim \mathcal{N}(0,1)$) of $(X_0,X_n)$.
Consider the strongly stationary Gaussian process $\{X_n\}_{n\in\N}$ where $X_n\sim \mathcal{N}(0,1)$ and the copula of $(X_r,X_s)$ is $C_{\rho_{|r-s|}}$ (such a construction is always possible for Gaussian processes, since all finite dimensional copulas can be taken $n$-dimensional Gaussian copulas with the appropriate covariance structure). Some examples are as follows.
\begin{itemize*}
\item For $q\in\Ns$, and $\{\vartheta_n\}_{n=0}^q$, the parameterization
\begin{equation*}\rho_n=\frac{1}{1+\vartheta_1^2}\sum_{j=0}^{q-|n|}\vartheta_j\vartheta_{j+|n|}\delta(|n|\leq q)\end{equation*}
determines a Gaussian MA$(q)$ process (cf. Brockwell and Davis, 1991, p.78 with the adequate adaptations).
\item For $|\varphi|<1$, the parameterization $\rho_n=\varphi^{|n|}$ determines a Gaussian AR(1) process (as noted in Joe, 1997).
\item The parameterization $\rho_n=2^{-n}(1+0.75n)$, for all $n\in\Ns$, determines a Gaussian ARMA$(2,1)$ process defined by $(1-\mathcal{B}+0.25\mathcal{B}^2)X_t=(1+\mathcal{B})Z_t$, $t\in\N$, for $\{Z_t\}_{t\in\N}$ i.i.d. $\mathcal{N}(0,32/3)$, where $\mathcal{B}$ denotes the backward shift operator (cf. Brockwell and Davis, 1991, p.92 with the adequate adaptations).
\item For $d\in(-0.5,0.5)$, the parameterization
\begin{equation}\label{arfpar}\rho_n=\prod_{k=1}^n\frac{k-1+d}{k-d}\sim\frac{\Gamma(1-d)}{\Gamma(d)}\,n^{2d-1},\end{equation}
for $n\in\Ns$, determines a Gaussian ARFIMA$(0,d,0)$ process (see Lopes, 2008), which is a long-range dependent process.
\item For $\{c_k\}_{k\in\N}$ a sequence of real numbers satisfying $\sum_{k\in\N}c_k^2<\infty$, the parameterization $\rho_n=\frac{\sum_{j\in\N}c_jc_{j+|n|}}{\sum_{k\in\N}c_k^2}$ determines a general Gaussian linear process with coefficients $\{c_k\}_{k\in\N}$.
\end{itemize*}
 \end{exe}
\begin{rmk}\label{cor-cov}
We notice that in Example \ref{gauss}, we assume that $X_n\sim\mathcal N(0,1)$. But in the usual definition of an ARFIMA$(0,d,0)$ process with standard normal innovations, for instance, the marginals are distributed as $\mathcal N\big(0,\sum_{k\in\N}c_k^2\big)$, where the sequence $\{c_k\}_{k\in\N}$ are the coefficients of the MA$(\infty)$ representation of the process. So that, in practice, one has to be careful to distinguish the meaning of the parameter $\rho_n$, which can be misleading.
\end{rmk}
The general framework of Theorem \ref{gen} can be extended to cover the case of parametric families of copulas for which the parameter space $\Theta\subseteq\R^k$, $k\in\Ns$.  Let $\{C_{\bs \t}\}_{\bs\t\in\Theta}$,
$\Theta\subseteq\R^k$, be a family of copulas for which $C_{\bs\t}$ is twice continuously differentiable with respect to $\bs\t$ on an open neighborhood $U\subseteq\Theta$ of a point $\bs a=(a_1,\cdots,a_k)^\tr\in\mathrm{int}(\Theta)$. Recall that the \emph{differential} of $C_{\bs\t}$ with respect to $\bs\t$ at $\bs a\in\R^k$ is the linear functional $d_{\bs\t}C_{\bs a}(u,v):\R^k\rightarrow\R$ whose value at a point $\bs b=(b_1,\cdots,b_k)^\tr\in\R^k$ is
\begin{equation*}d_{\bs\t}C_{\bs a}(u,v)\cdot\bs b=\sum_{i=1}^k\frac\d{\d\t_i}\,C_{\bs \t}(u,v)b_i\bigg|_{\bs\t=\bs a}.\end{equation*}
The \emph{second differential} of $C_{\bs\t}$ with respect to $\bs\t$ at $\bs a\in\R^k$ applied to $\bs b=(b_1,\cdots,b_k)^\tr\in\R^k$ is given by
\begin{equation*}d^2_{\bs \t}C_{\bs a}(u,v)\cdot\bs b^2=\sum_{i,j=1}^k \frac{\d^2}{\d\t_i\d\t_j}\,C_{\bs\t}(u,v)b_ib_j\bigg|_{\bs \t=\bs a}.\end{equation*}
With this formulation, in the following Theorem \ref{mp} we consider general decay of covariance in parametric families of copulas for which the space of parameter is a subset of $\R^k$, $k\geq 2$, allowing for some of the parameters to remain fixed.
\begin{thm}\label{mp}
Let $\{C_{\bs \t}\}_{\bs\t\in\Theta}$, for $\Theta\subseteq\R^{k+s}$ with non-empty interior, $k\in\Ns$ and $s\in\N$, be a family of parametric copulas for which there exists $\bs a\in\Theta'$ such that $\lim_{\bs \t\rightarrow\bs a}C_{\bs\t}(u,v)=uv$, for all $u,v\in I$. The limit is to be understood as the coordinatewise adequate lateral limits in case $a\notin\mathrm{int}(\Theta)$, also allowing for $s$ coordinates to remain fixed, that is, we allow for
\begin{equation*}\bs\t=(\t_1,\cdots,\t_k, \t^0_{k+1},\cdots,\t^0_{k+s})\longrightarrow(a_1,\cdots,a_k,\t^0_{k+1},\cdots,\t^0_{k+s})=\bs a.\end{equation*}
Also assume that there exists a set $D\subseteq\Theta$ with non-empty interior such that $\bs a\in D'$ and $C_{\bs\t}$ is twice continuously differentiable with respect to $\{\t_1,\cdots,\t_k\}$ in $D$. Let $\{F_n\}_{n\in\N}$ be a sequence of absolutely continuous distribution functions and define the sequences
\begin{equation*}K_1^{(i)}(n)=\iint_{I^2}\frac{1}{l_0(u)l_n(v)}\lim_{\bs\t\rightarrow \bs a}\frac{\d C_{\bs\t}(u,v)}{\d\t_i}\,dudv, \quad i=1,\cdots,k,\end{equation*}
and
\begin{equation*}K_2^{(i,j)}(n)=\iint_{I^2}\frac{1}{l_0(u)l_n(v)}\lim_{\bs\t\rightarrow \bs a}\frac{\d^2 C_{\bs\t}(u,v)}{\d\t_i\d\t_j}\,dudv, \quad  i,j=1,\cdots,k.\end{equation*}
Let $\{\bs\t_n\}_{n\in\Ns}$ be a sequence in $D$ converging to $\bs a$, with possibly $s$ fixed coordinates and let $\{X_n\}_{n\in\N}$ be a sequence of random variables such that $X_n\sim F_n$, $n\in\N$, and the copula associated with $(X_0,X_n)$ is $C_{\bs \t_n}$. Given a measurable function $R:\R\rightarrow\R$ satisfying $\lim_{n\rightarrow\infty}R(n)=0$, suppose that
\begin{equation}\label{mcond}
\sum_{i=1}^kK_1^{(i)}(n)(|\t_n^{(i)}-a_i|)\sim R(n)\quad\text{and}\quad \sum_{i,j=1}^kK_2^{(i,j)}(n)\big(\big|\t_n^{(i)}-a_i\big|\big) \big(\big|\t_n^{(j)}-a_j\big|\big)=o\big(R(n)\big).
\end{equation}
Then,  $\cov(X_0,X_n)\sim R(n)$ as $n$ goes to infinity.
\end{thm}
\begin{proof}
We present the proof for the case where $\bs a\notin\mathrm{int}(\Theta)$. The other cases are dealt analogously.  Let $\{\bs{\a}_m\}_{m\in\Ns}$ be an arbitrary sequence of parameters in $D$ such that $\bs{\a}_{m}\rightarrow \bs a$ (assuming the adequate lateral limit when necessary, allowing for $s$ coordinates to remain fixed). Applying Taylor's formula with Lagrange's remainder in $(\t_1,\cdots,\t_k)^\tr$ around $(a_1,\cdots,a_k)^\tr$ (the other parameters are fixed), there exists $\bs\t_0=(\t_0^{(1)},\cdots,\t_0^{(k)})^\tr$ such that $|\t_0^{(i)}-a_i|\leq|\t_i-a_i|$, for all $i=1,\cdots,k$, and
\begin{align}\label{assd}
C_{\bs\t}(u,v)&=\lim_{m\rightarrow\infty}\!C_{\bs\a_m}\!(u,v)+\!\lim_{m\rightarrow\infty}\!\!d_{\bs\t}C_{\bs\a_m}\!(u,v)(|\bs\t -\bs a|)+\frac{1}{2}\lim_{m\rightarrow\infty}\!\!d^2_{\bs\t}C_{\bs\a_m}\!(u,v)(\bs\t_0 -\bs a)^2\nonumber\\
&=uv+\sum_{i=1}^k\lim_{m\rightarrow\infty}\bigg[\frac{\d C_{\bs\t}(u,v)}{\d\t_i}\bigg|_{\bs\t=\bs\a_m}\bigg](|\t_i-a_i|)+\nonumber\\
& \hspace{2cm}+\frac{1}{2}\sum_{i,j=1}^k\lim_{m\rightarrow \infty}\bigg[\frac{\d^2C_{\bs\t}(u,v)}{\d\t_i\d\t_j}\bigg|_{\bs\t=\bs\a_m}\bigg](|\t_0^{(i)}-a_i|)(|\t_0^{(j)}-a_j|).
\end{align}
Substituting $\bs \t_n$ as in the enunciate in \eqref{assd}, we have
\begin{align}\label{asse}
C_{\bs\t_n}(u,v)&=uv+\sum_{i=1}^k\lim_{m\rightarrow\infty}\bigg[\frac{\d C_{\bs\t}(u,v)}{\d\t_i}\bigg|_{\bs\t=\bs\a_m}\bigg] (|\t_n^{(i)}-a_i|)+ \nonumber\\
& \hspace{2cm}+\frac{1}{2}\sum_{i,j=1}^k\lim_{m\rightarrow \infty}\!\bigg[\frac{\d^2 C_{\bs\t}(u,v)}{\d\t_i\d\t_j}\bigg|_{\bs\t=\bs\a_m} \bigg]\!\big(\big|\t_0^{(i)}(n)-a_i\big|\big)\!\big(\big|\t_0^{(j)}(n)-a_j\big|\big),
\end{align}
where $\t_0^{(i)}(n)\in\big(\min\{\t_n^{(i)},a_i\},\max\{\t_n^{(i)},a_i\}\big)\cup\{\t_n^{(i)}\}$, for $i=1,\cdots,k$. Let $\{X_n\}_{n\in\N}$ and $\{F_n\}_{n\in\N}$ be as in the enunciate. Hoeffding's lemma combined with \eqref{asse} yields
\begin{align*}
\cov(X_0,X_n)&=\sum_{i=1}^k\left[\iint_{I^2}\frac1{l_0(u)l_n(v)}\lim_{m\rightarrow\infty}\frac{\d C_{\bs\t}(u,v)}{\d\t_i}\bigg|_{\bs\t=\bs\a_m}dudv\right] (|\t_n^{(i)}-a_i|)+\nonumber\\
&\hspace{.8cm}+\frac{1}{2}\sum_{i,j=1}^k\bigg[\iint_{I^2}\frac1{l_0(u)l_n(v)}\lim_{m\rightarrow \infty}\frac{\d^2 C_{\bs\t}(u,v)}{\d\t_i\d\t_j}\bigg|_{\bs\t=\bs\a_m} \!\! dudv\bigg]\big(\big|\t_0^{(i)}(n)-a_i\big|\big)\big(\big|\t_0^{(j)}(n)-a_j\big|\big)\nonumber\\
&=\sum_{i=1}^kK_1^{(i)}(n)(|\t_n^{(i)}-a_i|)+\frac{1}{2}\sum_{i,j=1}^kK_2^{(i,j)}(n)\big(\big|\t_0^{(i)}(n)-a_i\big|\big) \big(\big|\t_0^{(j)}(n)-a_j\big|\big)\\
&\sim R(n)+o\big(R(n)\big)\sim R(n),
\end{align*}
by the hypothesis on $K_1^{(i)}$ and $K_2^{(i,j)}$.
\end{proof}
 Conditions \eqref{mcond} are general ones.  Observe that if there exists $i_0\in\{1,\cdots,k\}$ such that
 \begin{equation*}K_1^{(i_0)}(n)(|\t_n^{(i_0)}-a_{i_0}|)\sim R(n)\quad \text{and}\quad K_1^{(i)}(|\t_n^{(i)}-a_i|)=O\big(R(n)\big), \text{ for all }i\neq i_0,\end{equation*}
 and
 \begin{equation*}K_2^{(i,j)}(n)\big(\big|\t_n^{(i)}-a_i\big|\big) \big(\big|\t_n^{(j)}-a_j\big|\big)=o\big(R(n)\big), \text{ for all }i,j=1,\cdots,k,\ i\leq j,\end{equation*}
then the conclusion of Theorem \ref{mp} holds. These conditions, which imply \eqref{mcond}, are usually simpler to verify.

Recall that the Fr\`echet-Hoeffding lower bound copula and upper bound copulas are given respectively by, $W(u,v)=\max\{u+v-1,0\}$ and $M(u,v)=\min\{u,v\}$. The importance of the copulas $W$ and $M$ relies in two facts. First, for any copula $C$, $W(u,v)\leq C(u,v)\leq M(u,v)$, for all $u,v\in I$. Second, the copula related to the random vector  $\big(X,h(X)\big)$ is $W$ if, and only if, $h$ is decreasing and it is $M$ if, and only if, $h$ is increasing. A family of copulas for which the copulas $\Pi$, $W$ and $M$ are particular (or limiting) cases is said to be \emph{comprehensive}.

Comprehensiveness is a highly desirable property for a family of copulas because this means that the family can model a broad range of dependence structures. The Frank copula of Example  \ref{frexe} and the Gaussian copula of Example \ref{gauss} are examples of comprehensive families. The FGM family from Example \ref{fgm} is not, but due to the simple analytical form of the copulas in this family, it is widely employed in the literature in modeling, in testing association  and in studying efficiency of nonparametric procedures (cf. Nelsen, 2006, p.78). However, this family is arguably too restrictive for most applications as, for  instance, the Kendall's $\tau$ dependence coefficient for this family  ranges in $[-2/9,2/9]$. Since the space of all bidimensional copulas is a convex space, a simple solution for this problem is to consider a new copula obtained from the convex combination of the FGM copula and another copula(s) presenting the desired complementary characteristics.
 \begin{exe}
Let $\{C_\gamma\}_{\gamma\in[-1,1]}$ denote the FGM family of Example \ref{fgm} and let $\{C_\delta\}_{\delta\in[1,\infty)}$ be the Euclidean family of copulas with generator $\p_\delta(t)=(1-t)^\delta$. This family comprehend the copulas of the form
\begin{equation*}C_\delta(u,v)\vcentcolon=\max\big\{1-\big[(1-u)^\delta+(1-v)^\delta]^\frac{1}{\delta},0\big\},\quad\text{ for }\delta\geq1.\end{equation*}
Particular cases of this family are $C_1=W$ and $C_\infty=M$. Define a new three-parameter comprehensive family of copulas by setting
\begin{align}\label{3cop}
C_{\bs\t}(u,v)=\alpha C_\gamma(u,v)+(1-\alpha)C_{\delta}(u,v),
\end{align}
where $\bs\t\vcentcolon=(\gamma,\alpha,\delta)^\tr\in[-1,1]\times[0,1]\times[1,\infty)$. Notice that $C_{\bs\t}=\Pi$ when $\bs\t=(0,1,\delta)^\tr$, for all $\delta\in[1,\infty)$, so that, in the notation of Theorem \ref{mp}, $\bs a\vcentcolon=(0,1,\delta)^\tr$. In order to exemplify the use of Theorem \ref{mp}, we shall analyze the compatibility of the family \eqref{3cop} with long-range dependence. For simplicity, let us fix $\delta_0=1$ and consider the triple $\bs\t=(\gamma,\alpha,\delta_0)^\tr$. Consider the family of \emph{triangular distribution functions} in $[a,b]$, denoted by Tr$(a,b)$, whose distribution function and density are given, respectively, by
\begin{equation*}F(x;a,b)=\left(\frac{x-a}{b-a}\right)^2\quad \mbox{and} \quad f(x;a,b)=\frac{2(x-a)}{(b-a)^2},\end{equation*}
for all $x\in[a,b]$. For  two bounded sequences of real numbers $\{a_n\}_{n\in\N}$ and $\{b_n\}_{n\in\N}$, with $b_n>a_n$ for all $n\in\N$, let  $\{F_n\}_{n\in\N}$ be a sequence of distribution functions such that $F_n$ is distributed as Tr$(a_n,b_n)$, for each $n\in\N$. In this case $l_n(x)=\frac{2\sqrt x}{b_n-a_n}$. Let us denote $\d_\gamma C_{\bs\t}(u,v)\vcentcolon=\frac{\d C_{\bs\t}(u,v)}{\d\gamma}$ and similarly for the derivative with respect to $\alpha$.
The first derivative of $C_{\bs\t}$ with respect to $\gamma$ and $\alpha$ are
\begin{equation*}\d_\gamma C_{\bs\t}(u,v)=\alpha uv(1-u)(1-v)\quad\text{ and} \quad\d_\alpha C_{\bs\t}(u,v)=C_{\gamma}(u,v)-C_{\delta_0}(u,v), \quad\end{equation*}
so that
\begin{align*}
K_1^{(1)}(n)&=\frac14\left(\int_{I}\sqrt{u}(1-u)du\right)^{\!\!2}\!\!(b_0-a_0)(b_n-a_n)=\frac{4}{15^2}(b_0-a_0)(b_n-a_n)
\end{align*}
and
\begin{align*}
K_1^{(2)}(n)&=\frac{1}4\bigg[\left(\int_I\sqrt u\, du\right)^{\!\!2}-\iint_{I^2}\frac{W(u,v)}{\sqrt{uv}}\,dudv\bigg](b_0-a_0)(b_n-a_n)=k_0(b_0-a_0)(b_n-a_n),
\end{align*}
where $k_0=\frac{1}{9}-\frac{\pi}8$. As for the second derivative, since $C_{\bs\t}$ is a linear function of $(\gamma,\alpha)$,
\begin{equation*}\frac{\d^2 C_{\bs\t}(u,v)}{\d\gamma^2}=\frac{\d^2 C_{\bs\t}(u,v)}{\d\alpha^2}=0,\quad\text{and}\quad\frac{\d_{\gamma,\alpha}^2 C_{\bs\t}(u,v)}{\d\gamma\d\alpha}=uv\ln(u)\ln(v),\end{equation*}
so that
 \begin{equation*}K_2^{(1,2)}(n)=\frac19(b_0-a_0)(b_n-a_n)\qquad\text{and}\qquad K_2^{(1,1)}(n)=K_2^{(2,2)}(n)=0. \end{equation*}
 Therefore, if we choose $b_n-a_n\in\L$, with $b_n-a_n\nrightarrow0$, $\gamma_n\sim L(n)n^{-\beta}$ and $\a_n-1\sim o(n^{-\beta})$ (or vice-versa), for $\beta\in(0,1)$ and $L\in\L$, then $(C_{\bs\t},F_n)$ is compatible with long-range dependence. If we take $a_n=a_0$ and $b_n=b_0$, $a_0<b_0$, instead, we still obtain the result.
  \end{exe}
\begin{exe}
Recall that the \emph{Archimax} family of copulas (Cap\'era\`a et al., 2000; Mari and Kotz, 2001) comprehend copulas of the form
\begin{equation*}C_{\p,A}(u,v)=\p^{-1}\left(\big(\p(u)+\p(v)\big)A\left(\frac{\p(u)}{\p(u)+\p(v)}\right)\right),\end{equation*}
where $\p$ is an Archimedean generator and $A$ is a dependence function of an EVC family. Notice that when $\p\equiv1$ we obtain the Euclidean family with generator $\p$ and when $\p(t)=-\ln(t)$, we obtain the EVC family with dependence function $A$. Consider
\begin{equation*} \p_\t(t)=\ln\left(\frac{1-\t(1-t)}{t}\right) \ \ \ \mbox{ and } \ \ \ A_{\a,\beta}(t)=1-\min\{\a t,\beta(1-t)\},\end{equation*}
for $\theta\in[-1,1)$ and $\a,\beta\in I$. Notice that $\p_\t$ generates the AMH family from Example \ref{amh} while $A_{\a,\beta}$ generates the \emph{Marshall-Olkin} family of copulas (Nelsen 2006; Mari and Kotz, 2001), comprehending copulas of the form
\begin{equation*}C_{\a,\beta}(u,v)=\min\{uv^{1-\a},vu^{1-\beta}\}.\end{equation*}
Also notice that $A_{0,\beta}=A_{\a,0}\equiv1$. Fixed $\beta_0\in(0,1)$, the Archimax copula generated by $\p_\t$ and $A_{\a,\beta_0}$ can be written as
\begin{equation*}
C_{\p_{\t},A_{\a,\beta_0}}(u,v)=\begin{cases}
\frac{(\t-1)uv^{1-\a}}{\t uv^{1-\a}-[1-\t(1-u)][1-\t(1-v)]^{1-\a}}, & \mbox{ if }\ u\leq \frac{(\t-1)v^{\frac{\beta_0}{\a}}}{v^{\frac{\beta_0}{\a}}-[1-\t(1-v)]^{\frac{\beta_0}{\a}}};\\
\frac{(\t-1)vu^{1-\beta_0}}{\t uv^{1-\beta_0}-[1-\t(1-v)][1-\t(1-u)]^{1-\beta_0}e^{\beta_0}}, & \mbox{ if }\ u> \frac{(\t-1)v^{\frac{\beta_0}{\a}}}{v^{\frac{\beta_0}{\a}}-[1-\t(1-v)]^{\frac{\beta_0}{\a}}}.
\end{cases}
\end{equation*}
Hence, Theorem \ref{mp}, does not apply since there is no set where the copula is twice differentiable on $(\a,\t)$ for all $(u,v)\in I^2$. Notice that the support of the copula depends on both, the parameters $(\a,\t)$ and the argument $(u,v)$.
\end{exe}
\section{Conclusions and Final Remarks}
In this work, we study the problem of constructing stochastic processes with a predetermined decay of covariance by parameterizing a family of copulas and the process' marginals. Although the main interest in practice lies on stationary  processes, the theory proposed here covers both stationary and non-stationary processes and allow for arbitrary decay of covariance.  We show that the proposed methodology is free from the compatibility problem. The proof is novel and is based on the concept of monotonic regular indivisible index classes.

We present several examples to illustrate the theory, including the widely applied Euclidean, Extreme Value and Gaussian family of copulas. We show how the theory blends with some common applied time series models such as the large class of ARFIMA processes.
\subsubsection*{Acknowledgements}
G. Pumi was partially supported by CAPES/Fulbright Grant BEX 2910/06-3 and by CNPq-Brazil. S.R.C. Lopes research was partially supported by CNPq-Brazil, by CAPES-Brazil, by Pronex {\it Probabilidade e Processos Estoc\'asticos} - E-26/170.008/2008 -APQ1 and also by INCT {\it em Matem\'atica}.

\end{document}